\documentclass[11pt,fleqn,titlepage]{article}
\usepackage{latexsym,amsmath,amssymb,amsfonts,epsfig}
\usepackage[all]{xy}

\providecommand{\ra}{\rightarrow}

\providecommand{\CC}{{\mathbb{C}}}
\providecommand{\RR}{{\mathbb{R}}}

\providecommand{\ZZ}{{\mathbb{Z}}}

\providecommand{\QED}{{\hfill{} $\Box$}}

\providecommand{\DD}{{\cal D}} %% domain of operator
\providecommand{\GG}{{\cal G}}
\providecommand{\HH}{{\cal H}}
\providecommand{\KK}{{\cal K}}
\providecommand{\calS}{{\cal S}}

\providecommand{\THM}{{\mathbb{T}_H M}}

\providecommand{\THG}{{\mathbb{T}_H} G}

\providecommand{\Lg}{{\mathfrak g}}

\newcommand{\ang}[1]{\langle #1 \rangle} % <...>

\newtheorem{definition}{Definition}
\newtheorem{lemma}[definition]{Lemma}
\newtheorem{proposition}[definition]{Proposition}

\newtheorem{theorem}[definition]{Theorem}

\providecommand{\PSI}{{\overline{\Psi}}}
\providecommand{\EE}{{\cal E}}

\begin{document}

\title{The Atiyah-Singer Index Formula for Subelliptic Operators on Contact Manifolds, Part II}
\author{Erik van Erp \\ The University of Pennsylvania \\ erikerp@math.upenn.edu}
\date{April 2006}
\maketitle

%\tableofcontents

%\chapter{The index theorem for Fredholm operators in the Heisenberg calculus on contact manifolds}

\begin{center}{\Large \bf Abstract}\end{center}
\vskip 12pt

\noindent 
We present a new solution to the index problem for hypoelliptic operators in the Heisenberg calculus on contact manifolds,
by constructing the appropriate topological $K$-theory cocycle for such operators. Its Chern character gives a cohomology class to which the Atiyah-Singer index formula can be applied. 
Such a $K$-cocycle has already been constructed by Boutet de Monvel for Toeplitz operators, and, more recently, by Melrose and Epstein for the class of Hermite operators. Our construction applies to general hypoelliptic pseudodifferential operators in the Heisenberg calculus.

As in the Hermite Index Formula of Melrose and Epstein, our construction gives a vector bundle automorphism of the symmetric tensors of the contact  hyperplane bundle. This automorphism is constructed directly from the invertible Heisenberg symbol of the operator, and is easily computed in the case of differential operators.

\newpage

\section{Introduction}

In a previous paper (see [Er1]) we presented an index theorem for hypoelliptic differential operators on a compact contact manifold.
The hypoelliptic operators we considered were `elliptic' in the Heisenberg calculus.
Locally, every contact manifold is isomorphic to the Heisenberg group with its canonical contact structure.
As Stein and Folland showed [FS1], on the Heisenberg group a differential operator $P$ can be usefully approximated, at each point $m$, by a right-invariant differential operator $P_m$, obtained by `freezing coefficients'.
One can do the same on a contact manifold. 
The family of `model operators' $P_m$ thus obtained forms the symbol $\sigma_H(P)$ of $P$ in the Heisenberg calculus. 
If all model operators $P_m$ are hypoelliptic (which is easily verified by the Rockland condition) then $P$ is hypoelliptic, and, if $M$ is compact, $P$ is Fredholm.

In [Er1] we showed how to associate a $K$-theory element 
\[ [\sigma_H(P)]\in K^0(T^*M) \]
to a hypoelliptic Heisenberg symbol, and we proved that, with this choice of the $K$-theory symbol, the Atiyah-Singer formula computes the index of $P$,
\[ {\rm Index}\, P = \int_{T^*M} {\rm Ch}(\sigma_H(P)) \wedge{\rm Td}(M) .\] 
The construction of the $K$-cocycle $[\sigma_H(P)]$ proceeded in two steps.
First we obtained an element in the analytic $K$-theory of a noncommutative $C^*$-algebra $C^*(T_HM)$.
The groupoid $T_HM$ is just the tangent bundle $TM$, but each tangent space $T_mM$ carries a Heisenberg group structure encoding the contact structure of $M$.
The groups $G_m=T_mM$ are called {\em osculating groups} (see section \ref{section:two} below).
The $C^*$-algebra $C^*(T_HM)$ is the convolution algebra of this groupoid.
To obtain an element in $K^0(T^*M)$ we invoked the canonical isomorphisms
\[ K_0(C^*(T_HM)) \cong K_0(C^*(TM)) \cong K^0(T^*M) .\]
However, the first of these isomorphisms is highly nontrivial, and it is impossible to compute the element $[\sigma_H(P)]\in K^0(T^*M)$ in concrete examples in this way.

In the present paper we present the solution to this problem.
We compute the $K$-theory element $[\sigma_H(P)]$ as an explicit $K$-cocycle
\[ [\sigma_H(P)] = [a(P), V^N] \in K^1(M) .\]
Note that we work here in $K^1(M)$, instead of $K^0(T^*M)$. These groups are identical for contact manifolds that have a global contact form, which is a mild assumption we will make in this paper.
The vector bundle $V^N$ is associated to the hyperplane bundle $H$ (equipped with an almost complex structure),
\[ V^N = \bigoplus_{j=0}^N {\rm Sym}^j H^{1,0} ,\]
and the vector bundle automorphism $a(P)$ is constructed, at each point $m\in M$, from the model operators $P_m$, as the quotient
\[ a(P)_m = \pi_m(P_m)\pi_m^{op}(P_m)^{-1} .\]
Here $\pi_m$ denotes an explicit faithful irreducible representation of the osculating Heisenberg group $G_m$ on an infinite dimensional Hilbert space (the Bargmann-Fok space) that contains $V_m^N$.
The representations $\pi_m$ vary continuously along $M$,
and $\pi_m^{op}$ denotes the conjugate representation of $\pi_m$.

Recall that hypoellipticity of $P$ in the Heisenberg calculus is verified by the Rockland condition for the model operators $P_m$, which states  that  $\pi_m(P_m)$ and $\pi_m^{op}(P_m)$ are invertible. We will see that, as long as $N$ is sufficiently large, the quotient $a(P)_m$ is still invertible on the finite dimensional subspace $V_m^N$, and the $K$-theory class thus defined is independent of the choice of $N$.  

\vskip 6pt

It is interesting to compare our result with the index formulas found in the work of Melrose and Epstein [EM2], [Ep], [ME].
For a restricted class of operators, namely those for which 
\[ \pi_m(P_m^{op}) = 1 ,\]
our formula is equivalent to the `Hermite Index Formula' derived by Melrose and Epstein.
We merely recover this result by different methods.

However, in the solution of the general index problem, Melrose and Epstein introduce an `extended Chern form' associated with the Heisenberg symbol of a hypoelliptic operator. The construction of this Chern form is quite involved, and not very explicit. 
More problematically, it does not yield a closed de Rham form, i.e., the extended  Chern form is not really a `character'. To remedy this situation several extra terms appear in the final index formula in [EM2], which, as a result, becomes somewhat unwieldy.

We accomplish here for the general case what Melrose and Epstein had achieved for the `Hermite operators', which is an explicit construction of a cocycle in $K^1(M)$, whose Chern character can be computed by the usual Chern-Weil formalism, to which the classical Atiyah-Singer index formula applies. As a result we do not need any of the extra correction terms that appear in the formulas of Melrose and Epstein.
In particular, for {\em differential operators} this cocycle is easily computed, as we will show.

\vskip 6pt

Section \ref{section:main} contains the main results.
Sections \ref{section:two},  \ref{section:three}, and \ref{section:five} review the necessary background material.
Most of this material is taken from [Ep] and [EM2], which are the best sources for a detailed discussion of these facts.
A new result, interesting in its own right, is presented in section \ref{section:four}.
We prove the existence of a short exact sequence for the Heisenberg calculus,
\[ 0 \to \KK \to \PSI^0_H \to \overline{\calS}_H \to 0 .\]
Here $\PSI^0_H$ denotes the norm-closure of the order zero operators in the Heisenberg calculus (as bounded operators on $L^2(M)$),
and $\overline{\calS}_H$ is the closure, in an appropriate norm, of the algebra of principal symbols.

The three final sections give the proof of Theorem \ref{thm:main}. 
Here we make extensive use of the groupoids introduced in [Er1], but they are employed in a new way.
The tangent groupoid  method is merged here with the short exact sequence approach to index theory.

\vskip 6pt

We want to thank Charles Epstein for many valuable discussions.

\section{Contact manifolds and osculating groups}\label{section:two}

In this paper $M$ denotes a closed (compact, no boundary) contact manifold of dimension $2n+1$.
We let $H\subseteq TM$ denote the contact hyperplane bundle,
while $N$ denotes the quotient line bundle $N=TM/H$.
We will assume, for simplicity, that we can fix a {\em global} contact form $\theta$, with $\theta(H)=0$.
The choice of a global contact form $\theta$ corresponds to a section in the dual bundle $N^*$,
which, in turn, is equivalent to the choice of a trivialization $N^*\cong M\times \RR$.
Because $(d\theta)^n\theta$ is nondegenerate, the global contact form $\theta$ induces a symplectic form $\omega_m=d\theta_m$ on the fibers of $H$.

We define a group structure on the vector space $H_m\oplus N_m$, corresponding to the symplectic form $\omega_m$ on $H_m$.
Multiplication of elements in this group is given by
\[ (h,t)\cdot (h',t') = (h+h',t+t'+\frac{1}{2}\omega_m(h,h')) .\]
with $h,h'\in H_m, t,t'\in \RR$.
Here we identify $N_m=\RR$, corresponding to the trivialization of $N^*$.
The Lie algebra $\Lg_m=H_m\oplus N_m$ of this group is such that the exponential map is the identity.
While the symplectic structure on $H_m$ depends on the choice of $\theta$, the group structure on $G_m=H_m\oplus N_m$ is independent of this choice (simply verify what happens if we replace $\theta$ with $f\theta$).
 
The symplectic vector space $(H_m,\omega_m)$ is isomorphic to the standard symplectic space $\RR^{2n}$ with $\omega=\sum dx_idy_i$, where $(x,y) = (x_1,\ldots,x_n,y_1,\ldots,y_n)$.
Thus, the group $G_m$ is isomorphic to the Heisenberg group $G=\RR^{2n+1}$ with product
\[ (x,y,t)\cdot (x',y',t') = (x+x',y+y',t+t'+\frac{1}{2}\sum (x_iy'_i - x'_iy_i)) .\]
The groups $G_m$ are called the {\em osculating groups} of $M$.
The bundle of osculating groups forms a smooth groupoid, whose total space $H\oplus N$ is, of course, isomorphic to $TM$.
We denote this groupoid by $T_HM$.

The Lie algebras of the osculating groups are graded, with $H_m$ of degree $1$, and $N_m$ of degree $2$.
Correspondingly, these groups are equipped with a family of dilations $\delta_s$,
\[ \delta_s(h,t) = (sh,s^2t) .\]
The dilations $\delta_s$ are group automorphisms of $G_m$ (for $s>0$).

The notion of `homogeneity' of functions and distributions on the Heisenberg group conforms to the grading. 
For $\phi\in C_c^\infty(G)$, let $\phi_s(g) = \phi(\delta_s g)$.
Let $Q=2n+2$ denote the {\em homogeneous dimension} of $G$. 
We have $d(\delta_sg) = s^Qdg$, if $dg$ denotes Haar measure.
A distribution $k$ on $G$ is homogeneous of degree $d$ if $\ang{k,\phi_s} = s^{-Q-d}\ang{k,\phi}$.

\section{The Heisenberg calculus}\label{section:three}

In this section we briefly review the Heisenberg calculus.
For details see [Ep], [CGGP], [BG], [Ta1].

Like the classical pseudodifferential calculus, the Heisenberg calculus consists of linear operators 
\[ P\;\colon\; C^\infty(M)\to C^\infty(M) ,\]
whose Schwartz kernel $K\in \DD'(M\times M)$ is smooth off the diagonal in $M\times M$, and has an asymptotic expansion $K\sim \sum K^j$ near the diagonal.
The defining feature of the Heisenberg calculus is that the terms in the expansion of $K$ are homogeneous with respect to the `parabolic' dilations $\delta_s$ of the osculating groups.
Let's make this precise.

Because $P\colon C^\infty\to C^\infty$ we can identify the kernel $K$ with a smooth function 
\[ K\in C^\infty(M,\DD'(M)). \]
In other words, if we write  (formally), 
\[ P\phi(m) =  \int K(m,m')\phi(m') dm' = \ang{K_m,\phi} ,\]
then $K_m$ is a smooth family of distributions on $M$.
%i.e., the function $m\mapsto \ang{K_m,\phi_m}$ is smooth for every $\phi\in C^\infty(M\times M)$, where $\phi_m(m') = \phi(m,m')$. 

Now choose Darboux coordinates in an open set $U\subseteq M$,
i.e., identify $U$ with an open subset in the Heisenberg group $G$,
such that the contact form $\theta$ on $U$ is the pull-back of the canonical contact form on $G$.
Then let $k_m$ be the distribution on $G$ obtained by a translation of $K_m$ as follows \[ k_m(g) = K_m(g^{-1}m) .\]
Of course, $k_m$ will only be defined in a neighborhood of $0\in G$,
but this suffices to determine its asymptotic behaviour.
Because $K$ is smooth off the diagonal, each $k_m$ is regular, i.e., it is a smooth function when restricted to $G\setminus \{0\}$.
Now, locally, and up to a smoothing operator, we have the equality
\[ P\phi(m) = \int K(m,g)\phi(g)dg = \int k_m(mg^{-1})\phi(g)dg  = (k_m\ast \phi)(m) ,\]
where $\ast$ denotes convolution on $G$.
We see that the operator $P$ near the point $m$ can de approximated by a convolution operator on $G$. This was the basic idea introduced by Stein and Folland in [FS1].

We say that $P$ is a pseudodifferential operator of order $d$ in the Heisenberg calculus if
each distribution $k_m$ has an asymptotic expansion,
\[ k_m \sim k_m^0+k_m^1+k_m^2+\cdots ,\]
such that (1) each term $k_m^j$ in the expansion is regular and homogeneous of degree $d-j$, (2) each $k_m^j$ is a smooth family of distributions, while $k^j_m(g)$ is smooth in $(m,g)$ for $g\ne 0$ and (3) the remainder term
\[ k_m - \sum_{j=0}^N k_m^j = R_N \] 
becomes more and more smooth as $N$ grows (precisely: for every $l$ there exists $N$ such that $R_N\in C^l(G\times G)$).
This condition is independent of the choice of Darboux coordinates, and it suffices to verify the asymptotics in an open cover of $M$.
We denote the set of such operators by $\Psi_H^d(M)$.

The differential $T_0G\to T_mG$ of the map $G\to M\colon g\mapsto g^{-1}m$
induces an isomorphism of Lie algebras $\Lg\to \Lg_m=T_mG = H_m\oplus N_m$.
The exponential of this map gives an isomorphism of the osculating group $G_m$ at $m$ with $G$.
Thus, we can identify the principal part $k_m^0$ in the expansion of $k_m$ with a regular, homogeneous distribution on the osculating group $G_m$.
An easy calculation shows that the distribution $k_m^0$ is invariantly defined on $G_m$, independent of the choice of Darboux coordinates on $M$.

The right invariant, homogeneous operator $P_m$ on $G_m$ defined by convolution
\[ P_m\phi = k_m^0\ast \phi \]
is called the {\em `model operator'} of $P$ at $m$. 
The principal symbol $\sigma_H(P)$ of $P$ in the Heisenberg calculus is the smooth, right invariant family $\{P_m\}$ on the groupoid $T_HM$.  

Pseudodifferential operators in the Heisenberg calculus form a $\ast$-algebra.
If $P\in \Psi^a_H(M)$, $Q\in \Psi^b_H(M)$, then $PQ\in \Psi^{a+b}_H(M)$.
Also, each $P\in \Psi^a_H(M)$ has a formal adjoint $P^t\in \Psi^a_H(M)$
(with $\ang{P\phi,\psi} = \ang{\phi,P^t\psi}$ for  $\phi,\psi\in C^\infty$).
We have
\begin{align*}
  \sigma^{a+b}_H(PQ) & = \sigma^a_H(P)\ast \sigma^b_H(Q) ,\\
  \sigma^a_H(P^t) & = \sigma^a_H(P)^t .
\end{align*}  
For order zero operators the principal symbol map $\sigma_H$ is a $\ast$-homomorphism, and we have a short exact sequence
\[ 0\to \Psi^{-1}_H(M) \to \Psi^0_H(M) \stackrel{\sigma_H}{\longrightarrow} \calS_H \to 0 .\]
The algebra of symbols $\calS_H$ is the convolution algebra of smooth families of distributions on the osculating groups $G_m$ of $M$, that are regular and homogeneous of degree $-Q$ in each fiber.
In the analogous short exact sequence for classical pseudodifferential operators, we would have $\calS \cong C^\infty(S^*M)$. But in this case $\calS_H$ is a noncommutative algebra.

\section{A short exact sequence for the Heisenberg calculus}\label{section:four}

Our aim in this section is to complete the short exact sequence associated to the principal symbol map in the Heisenberg calculus to an exact sequence of $C^*$-algebras.

As in the classical calculus, the algebra $\Psi_H^0$ consists of bounded operators on $L^2(M)$, and we denote the norm closure by $\PSI_H^0$.
By the Rellich Lemma in the Heisenberg calculus, the norm closure of the ideal $\Psi^{-1}_H$ is the ideal of compact operators $\KK(L^2(M))$.

Let $\sigma = \{P_m,m\in M\}$ be a smooth family of homogeneous order zero model operators.  
We define a norm on the $\ast$-algebra $\calS_H$ by
\[ \|\sigma\|_\infty = \sup_{m\in M} \| P_m \|. \]
%\[ \|\sigma_H(P)\|_\infty = \sup_{m\in M, \pi\in \hat{G}_m} \| \pi(P_m) \|. \]
The right invariant model operators $P_m$ of an order zero operator $P$ are bounded as operators on $L^2(G_m)$, and $\|P_m\|$ denotes the operator norm.
The closure $\overline{\calS}_H$ of $\calS_H$ in this norm is a $C^*$-algebra.

The following proposition shows that the principal symbol map is continuous with respect to this norm.

\begin{proposition}
Let $P_m$ be the model operator of an order zero operator $P$ at $m\in M$.
Consider $P$ as a bounded operator on $L^2(M)$, and $P_m$ as a bounded operator on $L^2(G_m)$.
Then we have an inequality of operator norms,
\[ \|P\| \ge \|P_m\| .\]
\end{proposition}
{\bf Proof.}
We can think of $P_m$ and $P$ as operators on $G$, with $m=0$.
Let $V$ be a neighborhood of $0\in G$, and $h$ its characteristic function. Then,
\[ \|P\| \ge \|hPh\| \ge \|hP_mh\| - \|h(P-P_m)h\|.\]
The following two lemmas complete the proof.

\QED

\begin{lemma} If $P_m$ is an invariant order zero operator on $G$,
and $h$ the characteristic function of a neighborhood of $0\in G$, then
\[ \|hP_mh\| = \|P_m\|.\]
\end{lemma}
{\bf Proof.}
Because $P_m$ is of order zero,
we have $\|hP_mh\| = \|h_tP_mh_t\|$, where $h_t=h\circ \delta_t$.
As $t\to 0$, the support of $h_t$ blows up, and we have strong convergence $h_tP_mh_t \to P_m$.
This implies 
\[ \|hP_mh\| = \|h_tP_mh_t\|  \ge \|P_m\|. \]
The inverse inequality is trivial.

\QED

\begin{lemma}\label{lemma:key}
Given $\varepsilon>0$, there exists an open neighborhood $V\subseteq G$ of $0$ such that,
if $h$ denotes the characteristic function of $V$, then
\[ \|h(P-P_m)h\| \le \varepsilon .\] 
\end{lemma}
%%%% Comment: the editors of Annals of Math found a gap in the original proof of this lemma.
%%%% My original (incomplete) proof is preserved below.
%%%% This new proof was completed early February 2009.
{\bf Proof.}
Let $P_0$ denote the principal part of $P$. The definition of the operator $P_0$ only makes sense locally, and depends on a choice of coordinates.
In any case, $P-P_0$ is an order $-1$ operator whose kernel is an $L^1_{loc}$ function $a(x,y)$. 
Therefore we can arrange
\begin{align*} 
 \int |h(x) a(x,y) h(y)| dy  & \le \varepsilon ,\\
 \int |h(x) a(x,y) h(y)| dx  & \le \varepsilon ,
\end{align*} 
by choosing $V$ sufficiently small.
The first inequality holds for fixed $x\in G$, the second for fixed $y\in G$.
The combined inequalities imply the desired estimate for $P-P_0$,
\[ \|h(P-P_0)h\|\le \varepsilon.\]
A similar estimate holds for $P_0-P_m$, but for different reasons.
The difference $P_0-P_m$ is an order zero operator,
and its kernel 
\[ b(x,y) = k_x(y^{-1}) - k_m(xy^{-1})\]
is homogeneous of degree $-Q$ in the $y$-variable.
In [FS1], Folland and Stein give a proof of $L^2$-boundedness of singular integral operators
with such kernels, based on an application of the Cotlar-Stein lemma.
As is shown there, the norm of such an operator will be small if the kernel function $b(x,y)$ is uniformly small for all $y$-values on the unit spheres $|||y|||=1$.
(This fact is not stated explicitly, but implied by the details of the proof of Lemma 15.6 in [FS1]). 
Here $|||.|||$ denotes a homogeneous norm on $G$, i.e. a smooth function
for which $|||\delta_s y||| = s|||y|||$.
Since at $x=0$ we have $b(0,y)=0$---because, by definition of the model operator $P_m$, the kernels of $P_m$ and $P_0$ agree at $x=0$---it follows that $b(x,y)$ is close to $0$ for small $x$ and all $|||y|||=1$.
Therefore for sufficiently small $V$,
\[ \|h(P_0-P_m)h\|\le \varepsilon,\]
as desired.

\QED

%%%% The old, incomplete proof
%%{\bf Proof.}
%%Let $k(x,y)$ and $k_m(x,y)=k_m(xy^{-1})$ denote the Schwartz kernels of $P$ and $P_m$ respectively.
%%For small $xy^{-1}$, $x\ne y$, we have the inequality
%%\[ |k(x,y) - k_m(x,y)| \le C\,|||xy^{-1}|||^{-Q+1}  .\]
%%The triple bars $|||\cdot|||$ denote a `homogeneous norm' on $G$, i.e., a metric that satisfies $|||\delta_sz||| = s|||z|||$.
 %%Because $z\mapsto |||z|||^{-Q+1}$ is an $L^1_{loc}$ function,
%%we can arrange
%%\begin{align*} 
%% \int |h(x)(k(x,y)-k_m(x,y)) h(y)| dy  & \le \varepsilon ,\\
%% \int |h(x)(k(x,y)-k_m(x,y)) h(y)| dx  & \le \varepsilon ,
%%\end{align*} 
%%by choosing $V$ sufficiently small.
%%The first inequality holds for fixed $x\in G$, the second for fixed $y\in G$.
%%The combined inequalities imply the required estimate for the operator norm.
%%
%%\QED

We obtain a short exact sequence of $C^*$-algebras 
\[ 0 \to \KK \to \PSI^0_H \to \overline{\calS}_H \to 0 ,\]
with the usual corollaries:
(1) an order zero operator $P\in \PSI^0_H$ is Fredholm if and only if its symbol $\sigma_H(P)$ is invertible in the algebra $\overline{\calS}_H$, (2) an invertible symbol determines a class in the $K$-theory group 
\[ [\sigma_H(P)] \in K_1(\overline{\calS}_H) ,\]
and (3) the boundary map in $K$-theory
\[ \partial \;\colon\; K_1(\overline{\calS}_H) \to K_0(\KK) = \ZZ \]
sends the symbol of $P$ to the Fredholm index of $P$,
\[ {\rm Index}(P) = \partial[\sigma_H(P)] .\]
It can been shown that if each model operator $P_m$ is invertible as a bounded operator on $L^2(G_m)$, with two-sided inverse $Q_m$, then the family $\{Q_m\}$ determines an element in $\calS_H$, and hence $\sigma_H(P)$ is invertible.
%Moreover, this is equivalent to the {\em Rockland condition}, which states that for each nontrivial irreducible unitary representation $\pi$ of $G_m$ the operator $\pi(P_m)$ is invertible.
([CGGP], Theorems 2.5.d and 8.1).

\section{Harmonic analysis on the osculating groups}\label{section:five}

The structures outlined in this section are described in detail in the work of Melrose and Epstein. We refer to [Ep] for the details.

Let $V\cong \RR^{2n}$ be a symplectic vector space with symplectic form $\omega$. 
The group $G=V\oplus \RR$ with multiplication
\[ (v,t)\cdot (v',t') = (v+v',t+t'+\frac{1}{2}\omega(v,v'))  \]
is isomorphic to the Heisenberg group.
(We are, of course, interested in the osculating groups $G_m$, with $V=H_m$, $\omega=d\theta_m$.)
Let $J\colon V\to V$ be a complex structure ($J^2=-1$) on $V$ `adapted' to the symplectic structure.
This means that $\omega(Ju,Jv) = \omega(u,v)$, and $\omega(Jv,v)>0$ if $v\ne 0$.

The complexified space $V\otimes \CC$ splits naturally into the $\pm \sqrt{-1}$ eigenspaces for $J$, 
\[ V\otimes \CC = V^{1,0} \oplus V^{0,1}. \]
Correspondingly, the complexification of the Lie algebra $\Lg$ of $G$ splits as a direct sum 
\[ \Lg \otimes \CC = V^{1,0}\oplus V^{0,1}\oplus \CC. \]
If we extend the bilinear form $\omega$ to a complex bilinear form on $V\otimes \CC$,
then the expression
\[ \ang{z,w} = 2i\omega(z,\bar{w}) \]
defines a hermitian inner product on $V^{1,0}$. 
Let $\{Z_1,\ldots,Z_n\}$ be a basis for $V^{1,0}$ that is orthonormal with respect to this inner product. 
The Lie bracket on $\Lg\otimes \CC$ is then given by
\[ [Z_j, Z_k] = [\bar{Z}_j, \bar{Z}_k] = 0 \;,\; [Z_j,\bar{Z}_k] =  \frac{1}{2i} \delta_{jk} .\]
The {\em Bargmann-Fok} space $\HH^{BF}$ is the Hilbert space of holomorphic functions on $V^{0,1}$ with inner product
\[ \ang{f,g} = \int f(z)\overline{g(z)} \,e^{-|z|^2} dz .\]
The Bargmann-Fok representation $\pi$ of the complexified Lie algebra $\Lg\otimes\CC$ on the space $\HH^{BF}$ is given by
\[ \pi(Z_j) = iz_j \;,\; \pi(\bar{Z}_j) = -i\frac{\partial}{\partial z_j} \;,\; \pi(1) = \frac{1}{2}i .\]
This induces an irreducible unitary representation $\pi$ of $G$.

Now consider the {\em anti-automorphism} of $G$ given by 
\[ {\rm op}\;\colon\; V\oplus \RR\to V\oplus \RR \;\colon\; (v,t) \mapsto (v,-t) .\]
For $\phi\in C_c^\infty(G)$ we let $\phi^{op}(v,t) = \phi(v,-t)$,
and for a distribution $k$ on $G$ we have $\ang{k^{op},\phi} = \ang{k,\phi^{op}}$.
If $k$ is homogeneous of degree $-Q$, then so is $k^{op}$, and $\widehat{k^{op}} = \hat{k}^{op}$.
If $k$ is invertible (for the convolution product) then so is $k^{op}$, because `op' is an anti-automorphism. 

In the remainder of this section we will prove the following result, which will be the key to our construction of a  $K$-cocycle from the model operators.
\begin{proposition}\label{prop:weyl}
Let $k$ be a regular distribution on the Heisenberg group $G=V\oplus \RR$ that is homogenous of degree $-Q$. 
Assume that $k$ has a two-sided inverse for the convolution product.
If $\pi$ denotes the Bargmann-Fok representation of $G$, then the operator
\[ \pi(k)\pi(k^{op})^{-1} - 1 \]
is compact. 
\end{proposition}
{\bf Proof.}
The proof of this proposition is an application of the {\em Weyl calculus}, which is intimately connected with the representation theory of the Heisenberg group (see [Ep]).

Let $a_+$ denote the restriction of a function $a$ on $V^*\times \RR^*$ to the hyperplane $V^*_+ = V^*\times \{1\}$,  i.e., $a_+(\xi)= a(\xi,1)$, $\xi\in V^*$.
A straightforward calculation shows that for two $L^1$ functions $f,g$ on $G$, the restriction of the Fourier transform of their convolution product $\widehat{f\ast g}$ to $V_+^*$ only depends on the restrictions of the Fourier transforms $\hat{f},\hat{g}$ to $V_+^*$, as follows,
\[ \widehat{(f\ast g)}_+(\xi) = (\hat{f}_+ \,\#\, \hat{g}_+)(\xi)
 = \int e^{i\omega(u,v)} \hat{f}_+(\xi+u)\hat{g}_+(\xi+v)dudv .\]
This sharp product $\#$ is the product of symbols in the Weyl calculus.
In this calculus, a smooth function $a$ on $V^*_{+}$ is quantized as a Hilbert space operator $q_W(a)$ (for an explicit formula, see [Ep]).
If $k$ is a regular distribution on $G$, homogeneous of degree $-Q=-2n-2$,
then $\hat{k}$ is regular and homogeneous of degree $0$ on $V^*\oplus \RR^*$, i.e.,
it is a smooth function on $G\setminus \{0\}$ that is constant on the parabolic rays $(s\xi,s^2\tau)$.
Therefore $\hat{k}_+$ is a smooth function  that extends continuously to the radial boundary $\partial V_+^* \cong S^{2n-1}$.
Thus,  $q_W(\hat{k}_+)$ is an order zero operator in the Weyl calculus.
In particular, $q_W(\hat{k}_+)$ is bounded.
The restriction of the Weyl symbol $\hat{k}_+$ to the boundary $S^{2n-1}$ of $V^*_+$ is the {\em principal} Weyl symbol $\sigma^0_W(\hat{k}_+)$. It satisfies
\[ \sigma^0_W(a\,\# b) = \sigma^0_W(a)\sigma^0_W(b) .\]
Because $\widehat{k^{op}} = \hat{k}^{op}$, while $\hat{k}$ is constant on the parabolic rays $(s\xi,s^2\tau)$ with $(\xi,\tau)\in V\oplus \RR$, one easily verifies that
\[ \sigma^0_W(\hat{k}_+)=\sigma^0_W(\hat{k}_+^{op}). \] 
In particular, if $k$ is invertible, then
$\sigma^0_W(\hat{k}_+)\sigma^0_W(\hat{k}_+^{op})^{-1} = 1$
implies that the operator 
\[ q_W(\hat{k}_+)q_W(\hat{k}_+^{op})^{-1} - 1 \]
is of order $-1$ in the Weyl calculus, and is therefore {\em compact}.

The hyperplane $V_+^*$ is one of the coadjoint orbits of the Heisenberg group $G=V\oplus \RR$.
For nilpotent groups, Kirillov theory establishes a one-to-one correspondence between irreducible unitary representations and coadjoint orbits.
The Bargmann-Fok representation corresponds to $V_+^*$, and we can, in fact, identify
\[ \pi(k) = q_W(\hat{k}_+) .\]

\QED

\section{The main result}\label{section:main}

In this section we construct a topological $K$-cocycle from the family $\{P_m\}$.
Once this is done, we state our main result, Theorem \ref{thm:main}, and give two applications. The remainder of this paper will be concerned with the proof of Theorem \ref{thm:main}.

The osculating groups $G_m$ are Heisenberg groups.
We fix a {\em global} contact form $\theta$ on $M$ (we assume this can be done),
and we identify $N^*=M\times \RR$ correspondingly.
Thus we can identify $G_m = H_m\oplus \RR$, where $H_m$ is a symplectic space with $\omega_m=d\theta_m$, and we can use the structures developed in the previous section. 
There exists an almost complex structure $J$ on $H$ adapted to the symplectic structure $\omega_m=d\theta_m$ in each of its fibers (see [Ep]).
At each point $m\in M$, we now have the Bargmann-Fok representation $\pi_m$ of $G_m$ on the Hilbert space of holomorphic functions on $H_m^{0,1}$, with Gaussian measure. We denote this Hilbert space by $V_m^{BF}$.

Let $\Gamma$ be the set of continous functions $f$ on the space $H^{0,1}$ for which the restriction to any of the fibers $H_m^{0,1}$ is a complex polynomial $f_m\in V^{BF}_m$. 
Clearly, since the family $f_m$ consists simply of polynomials with continuous coeficients, the Bargmann-Fok norms 
\[ \|f_m\|^2 = \int |f_m(z)|^2 \,e^{-|z|^2} dz \]
are continuous in $m$, so that $\Gamma$ determines a {\em continuous structure} on the family of Hilbert spaces $V_m^{BF}$.   
We denote the continuous field of Hilbert spaces $\{V_m^{BF}, m\in M\}$ with the continuous structure defined by the set $\Gamma$ by $V^{BF}$. 
Observe that the elements in $\Gamma$ can be identified with the sections in the finite dimensional vector bundles
\[ V^N = \bigoplus_{j=0}^N {\rm Sym}^j H^{1,0}.\] 
The fiber $V^N_m\subset V^{BF}_m$ corresponds precisely to the space of complex polynomials on $H^{0,1}_m$ of degree $\le N$.

We now build a $K^1$-cocycle from the invertible model operators $P_m$ of a Fredholm operator $P\in \Psi^0_H$. We need a preliminary result.

\begin{proposition}\label{prop:cont}
If $P_m$ is the family of model operators of an order zero operator $P$,
then $\{\pi_m(P_m)\}$ is a norm continuous family of operators on the continuous field $V^{BF}$.
\end{proposition}
{\bf Proof.}
We choose a local trivialization $G_m\cong G$ of the bundle of osculating groups,
which is a trivialization $H_m=V$ such that $\omega_m$ is a constant form on $V$.
The model operators $P_m$ are then identified with a family of operators on $L^2(G)$,
and the Bargman-Fok representation $\pi_m$ is the same at each point.

We have $k^0_m = \tilde{k}^0_m + a(m)\delta$, where $\tilde{k}^0_m$ is a principal value distribution on $G$, $\delta$ is the Dirac delta, and $a\in C^\infty(M)$. By assumption, $\tilde{k}^0_m(g)$ is smooth in $(m,g)$, for $g\ne 0$. 
We can estimate the operator norm $\|P_m-P_{m'}\|$ (for $m$ close to $m'$) in terms of the supremum of the function $k^0_m-k^0_{m'}$, restricted to the unit sphere.
(This is how one proves that a principal value distribution defines a bounded operator on $L^2(G)$, see [FS1].) 
We then see that the family $P_m$ itself, as operators on $L^2(G)$, is norm continuous

Because $\|P_m-P_{m'}\|\ge \|\pi(P_m-P_{m'})\|$,
it follows that the family $\pi(P_m)$ is norm continous as well.

\QED

Now let $P$ be an operator that has an invertible symbol in the Heisenberg calculus.
Such an operator is hypoelliptic, and Fredholm if $M$ is compact.
Consider the family $a(P) = \{a(P)_m,m\in M\}$ of invertible operators on $V^{BF}_m$, defined as
\[ a(P)_m = \pi_m(P_m)\pi_m(P_m^{op})^{-1} .\]
By Proposition \ref{prop:weyl} we have
\[ a(P)_m -1 \in \KK(V^{BF}_m) ,\] 
while Proposition \ref{prop:cont} implies that
\[ a(P) \in \KK^+(V^{BF}) .\]
Here $\KK(V^{BF}_m)$ denotes the $C^*$-algebra of compact operators on the Hilbert space $V_m^{BF}$, while $\KK^+(V^{BF})$ denotes the unitalization of the $C^*$-algebra of continuous sections in the field of $C^*$-algebras $\{\KK(V_m^{BF})\}$.
Note that $\KK(V^{BF})$ are the `compact operators', in the sense of Kasparov, on the Hilbert module $V^{BF}$. 
Therefore, since $a(P)$ is invertible, it defines an analytic $K^1$-cocycle
\[ [a(P),V^{BF}] \in K_1(C(M)) \cong K^1(M).\]
To obtain a topological cocycle, we `compress' $a(P)$ to a finite-dimensional vector bundle $V^N\subset V^{BF}$ (i.e., if $e_N$ denotes the projection of $V^{BF}_m$ onto $V^N_m$,  take $e_N a(P)_m e_N$). 
If $N$ is sufficiently large, $a(P)$ restricts to an automorphism of the vector-bundle $V^N$, and the cocylcle
\[ [a(P),V^N] \in K^1(M) \]
is equivalent to $[a(P),V^{BF}]$.

We can now state our main result.

\begin{theorem}\label{thm:main}
Let $P$ be an order zero pseudodifferential operator in the Heisenberg calculus
with invertible model operators.
Then
\[ {\rm Index}\, P = \lim_{N\to \infty} \int_M {\rm Ch}([a(P), V^N]) \wedge{\rm Td}(M) ,\]
where 
\[ {\rm Ch}\;\colon\; K^1(M)\to H^{\rm odd}(M) \]
denotes the classical odd Chern character.
The limit stabilizes for values of $N$ that are sufficiently large.
\end{theorem}
\noindent{\bf Remark.}  Restriction to a finite dimensional vector bundle $V^N$ is a triviality in the context of analytic $K$-theory.
It allows us to compute the Chern character by the classical means.
It is precisely at this point that our approach departs from that of Melrose and Epstein.
Melrose and Epstein attempt to define a Chern form for each of the families $\pi_m(P_m)$ and $\pi_m(P_m^{op})$ by itself (equivalent to their $\sigma^H(P)(+1)$ and $\sigma^H(P)(-1)$), i.e., before taking their quotient, as we do
(in the context of their work their is no reason to think that one {\em should} take a quotient). Since neither of these elements by itself defines a $K^1$-cocycle, one can see why such an attempt must meet with unavoidable difficulty.

For Hermite operators, i.e., operators for which $\pi_m(P^{op}_m) = 1$ (see section \ref{section:nine} below),
there is no need to pass to the quotient $a(P)$, and the methods of [EM2] are sufficient. Their {\em Hermite Index Formula} is, essentially, equivalent to our result. 
Our formula is valid in the general case.
\vskip 6pt

Before we turn to the proof of Theorem \ref{thm:main}
we show how it solves the index problem for hypoelliptic {\em differential} operators on a contact manifold.
Suppose $P$ is a pseudodifferential Heisenberg operator of order $d$ with Rockland model operators.
The choice of a complex structure in the fibers of $H$ fixes a sublaplacian $\Delta = -\sum (X_j^2+Y_j^2)$ on $M$ (up to lower order terms).
Then $P(1+\Delta)^{-d/2}$ is a Fredholm pseudodifferential operator of order zero in the Heisenberg calculus, with the same index as $P$.  

\begin{proposition}\label{prop:diffop}
Let $P$ be a differential operator of order $d$ in the Heisenberg calculus,
with Rockland model operators.
The family of operators
\[ a(P)_m= \pi_m(P_m)\pi_m(P_m^{op})^{-1} \]
defines  an invertible element in $\KK^+(V^{BF})$,
and we have
\[ [a(P)] = [a(P(1+\Delta)^{-d/2})] \in K^1(M)  .\]
Therefore, Theorem \ref{thm:main} holds for differential operators.
\end{proposition}
{\bf Proof.}
The Bargmann-Fok representation $\pi_m(\Delta_m)$ of $\Delta_m$  is the harmonic oscillator $Q_m$ on the Hilbert space $V^{BF}_m$.
The summands ${\rm Sym}^k H^{1,0}_m\subseteq V^{BF}_m$ are the eigenspaces for $Q_m$,
with eigenvalues $k+\frac{1}{2}n$ (see [Ep]).

For a general differential operator $P$,
the representation $\pi_m(P_m)$ on the space $V^{BF}_m$ of holomorphic functions on $H^{0,1}_m$  is a polynomial in $z_j$ and $\partial/\partial z_j$.
Therefore, it will map the eigenspace ${\rm Sym}^k H^{1,0}_m$  to a subspace of the finite sum
\[ \bigoplus_{|l|\le d} {\rm Sym}^{k+l} H^{1,0}_m .\]
In other words, $\pi_m(P_m)$ has `finite propagation' on the spectrum of $Q_m$.

Now, because $((P_m\Delta_m^{-d/2})^{op})^{-1} = (\Delta_m^{-d/2}P_m^{op})^{-1} =
(P_m^{op})^{-1}\Delta_m^{d/2}$, we find 
\[ a(P(1+\Delta^{-d/2}))_m =  \pi_m(P_m)Q_m^{-d/2}\pi_m(P_m^{op})^{-1}Q_m^{d/2}  .\]
Our spectral analysis of $\pi_m(P_m)$ shows that,
if we restrict to a finite dimensional  bundle $V^N$ for large $N$,
we can replace 
$Q_m^{-d/2}\pi_m(P_m^{op})^{-1}Q_m^{d/2}$ with $\pi_m(P_m^{op})^{-1}$
and get the same $K^1$-cocycle.

\QED

\noindent
{\bf Example: second order differential operators.}
Consider the second order differential operator
\[ P = \sum (X_j^2+Y_j^2) + i\beta T = \sum Z_j\bar{Z}_j - i(n-\beta) T .\]
Here $\beta$ is a complex-valued function. 
Melrose and Epstein derived an explicit index formula for these `twisted sub-Laplacians' as a corollary of their Hermite Index Formula (see Chapter 11 in [EM2]). 
We show here how their formula can be derived from Theorem \ref{thm:main}.

The model operators of $P$ are
\[ P_m      = \sum Z_j\bar{Z}_j - i(n-\beta(m))T ,\]
%\[ P_m^{op} = \sum \alpha_j(m)Z_j\bar{Z}_j - i(\sum \alpha_j(m) + \beta(m))T \]
and we find
\[ \pi_m(P_m)      = \sum z_j\frac{\partial}{\partial z_j} + \frac{1}{2}(n - \beta(m)) .\]
%\[ \pi_m(P_m^{op}) = \sum \alpha_j(m)z_j\frac{\partial}{\partial z_j} + \frac{1}{2}(\sum \alpha_j(m) + \beta(m)) .\]
The action of $\pi_m(P_m)$ on the Bargmann-Fok space is given by
\[ \pi_m(P_m) z^\alpha = (|\alpha| + \frac{1}{2}(n-\beta(m))\, z^\alpha .\]
%\[ \pi_m(P_m^{op}) z^\alpha = (|\alpha| + \frac{1}{2}\beta(m) +n))\, z^\alpha \]
To find $\pi_m(P_m^{op})$, simply replace $\beta$ with $-\beta$.
We see that $P$ is subelliptic (and Fredholm) if and only if $\beta$ does not take values in the exceptional set
\[ \Lambda = \{ \ldots, -n-4, -n-2, -n, n, n+2, n+4, \ldots \},\]
and by homotopy invariance of the Fredholm index, the index of $P$ must depend on the homotopy type of the map 
\[ \beta\;\colon\; M\to \CC\setminus \Lambda .\]
On the summands ${\rm Sym^k}H^{1,0}$ in $V^{BF}$, we have the scalar action 
\[ a(P) = \frac{n+2k-\beta}{n+2k+\beta} .\]
Clearly, we can ignore all values of $k$ for which both exceptional values $\pm(n+2k)$ lie in the unbounded component of the complement of the image $\beta[M]\subseteq \CC$.
In that case the automorphism $a(P)$ is homotopically trivial on ${\rm Sym^k}H^{1,0}$, and does not contribute to the Chern character of $[a(P),V^{BF}]$.
The formula in Theorem \ref{thm:main} reduces to a finite sum
\[ {\rm Index}\, P = \sum_k \int_M {\rm Ch}\left(\frac{n+2k-\beta}{n+2k+\beta}\right)\wedge {\rm Ch}({\rm Sym}^k H^{1,0}) \wedge {\rm Td}(M) ,\]
which is equivalent to the formula given in [EM2].

Of course, for general differential operators, the computation is not quite so easy, but it is the same in principle. 
By choosing a local representation of the manifold by Darboux coordinates, one represents $P$ as a sum of monomials of vector fields $Z_j,\bar{Z}_j,T$.
One easily finds the matrix entries for $\pi_m(P_m)$ and $\pi_m(P_m^{op})$ (by choosing the basis $z^\alpha$ of $V^{BF}$), which have finite propagation on the summands ${\rm Sym}^k H^{1,0}$ (see the proof of proposition \ref{prop:diffop}).
The entries of $\pi_m(P_m)\pi_m(P_m^{op})^{-1}$ (computed by means of elementary matrix algebra) converge to $1$ as $k$ gets large,
and so one can restrict the matrices to a finite dimension, yielding an explicit element in $K^1(M)$.

This may not sound like something that one would like to do in practice, but then again, the same is true for the calculation of the Chern character for elliptic operators!

\vskip 6pt
\noindent{\bf Example: Toeplitz operators.}
Melrose and Epstein take Boutet de Monvel's index theorem for Toeplitz operators as the starting point of their work, and work their way up from there to a general formula.
By contrast, in our case the Toeplitz index theorem appears as a corollary of the general formula.
We simply use some facts and tricks from [EM2], but in the reverse direction (from general theorem to specific case).

Let $M=\partial X$ be the smooth boundary of a strictly pseudoconvex domain $X\subset \CC^{n+1}$, with its natural contact structure.
The Hardy space $H^2(M)$ is the closure in $L^2(M)$ of smooth functions on $M$ that extend holomorphically to $X$,
and the Szeg\"o projector $S$ is the orthogonal projection of $L^2(M)$ onto $H^2(M)$.

The projector $S$ is not a pseudodifferential operator in the ordinary calculus,
but it {\em is} an order zero operator in the Heisenberg calculus.
The Bargmann-Fok representations $\pi_m(S_m)$ of the model operators $S_m$ of $S$
are rank one projections $s_m$.
In fact, $s_m$ is the projection onto the span of the vacuum vector (with eigenvalue $n/2$) of the harmonic oscillator $Q_m$,
which is the one-dimensional space $V^0_m\cong \CC$.
Moreover,
\[ \pi_m(S_m^{op})=0 .\]
(See [Ep], [EM2] for details.)

For a smooth function $a$ on $M$, the Toeplitz operator $T_a\colon H^2(M)\to H^2(M)$ is defined by 
\[ T_a = SaS. \]
For the purpose of index theory, we may replace the operator $T_a$ on $H^2(M)$ with 
the operator $P=SaS+(1-S)$ on $L^2(M)$.
We find $\pi_m(P_m) = a(m)s_m + (1-s_m)$, $\pi_m(P_m^{op}) = 1$, and therefore, 
\[ a(P)_m = a(m)s_m + (1-s_m) .\]
Because $a(P)_m$ is just the identity operator on the orthogonal complement of $V^0$ in $V^{BF}$, we can restrict $a(P)$ to the trivial line bundle
\[ V^0 = {\rm Sym^0} H^{1,0} = M\times \CC,\]
which leaves us with 
\[ [a(P),V^{BF}] = [a(P), V^0] = [a] .\]
Theorem \ref{thm:main} then reduces to the formula of Boutet de Monvel [Bo],
\[ {\rm Index}\, T_a = {\rm Index}\, P = \int_M {\rm Ch}(a) \wedge Td(M) .\]
{\bf Remark.} To get a meaningful index theorem for Toeplitz operators, we must consider {\em matrices} of such operators, or, alternatively, operators acting in sections of a trivial bundle $M\times \CC^N$.
(As stated here, we have always ${\rm Index}\,T_a = 0$.)
However, it is a trivial matter to formulate and prove Theorem \ref{thm:main} for such operators.

\section{A calculus on the parabolic tangent groupoid}

We now turn to the proof of Theorem \ref{thm:main}.
We start, in this section, by  extending the Heisenberg calculus on $M$ to a calculus of operators on the {\em parabolic tangent groupoid} $\GG = \THM$ described in [Er1].
We define our calculus in terms of asymptotic expansions of operator kernels. 
This approach to the Heisenberg calculus is developed in detail in a paper by Christ, Geller, Glowacki, and Polin [CGGP].
For proofs of the facts stated in this section we refer to [CGGP].
The results obtained there for the Heisenberg calculus generalize easily to the calculus we introduce here.

Recall that the parabolic tangent groupoid $\GG=\THM$ is a union of $T_HM$ with a family of pair groupoids $M\times M$,
\[ \GG = T_HM \cup M\times M\times (0,1] .\]
The family of pair groupoids $M\times M\times (0,1]$ is glued to the bundle $T_HM$ of osculating groups by blowing up the diagonal in $M\times M$ using the parabolic dilations of the Heisenberg group.
This is done locally, by identifying open subsets $U\subseteq M$ with open subsets in $G$.
The parabolic tangent groupoid $\THG$ of the Heisenberg group itself can be identified with the transformation groupoid $B\rtimes_\alpha G$, where $B = G\times [0,1]$, and $\alpha(g)(p,s) = (\delta_s(g)p,s)$, $g,p\in G$, $s\in [0,1]$.
(See [Er1] for details.)

Let
\[ K\in C_c^\infty(B,\EE'(G)) \]
be a compactly supported smooth family $K_x,x\in B$, of compactly supported distributions $K_x\in \EE'(G)$. 
We think of $K$ as a kernel on the groupoid $\GG=B\rtimes_\alpha G$.
We say that $K$ is {\em regular} if $K_x(g)$ is a smooth function when restricted to $\GG\setminus \GG^{(0)}$.

Using the triple notation $(x,g,y)$ for elements in $B\rtimes_\alpha G$
(where $x,y\in B$ denote target and source element, respectively, and $g\in G$),
convolution of two regular kernels $K,L$ is defined formally by
\[ (K\ast L)(x,g,y) = \int_G K(x,h,z)L(z,h^{-1}g,y)dh ,\]
or, equivalently,
\[ (K\ast L)_x (g) = \int_G K_x(h)L_{\alpha(h^{-1})x}(h^{-1}g)dh .\]
More precisely, for $\phi\in C^\infty(G)$,
\[ \ang{(K\ast L)_x,\phi} = \ang{K_x(h), \,\ang{ L_{\alpha(h^{-1})x}(h^{-1}g), \phi(g)}} .\]

\begin{proposition}
If $K,L$ are regular kernels on $\GG$, then so is $K\ast L$.
\end{proposition}
{\bf Proof.}
The proof is similar to that for convolution of regular kernels on a group.

First of all, for fixed $\phi\in C^\infty$ the expression $\ang{ L_{\alpha(h^{-1})x}(h^{-1}g), \phi(g)}$ is a smooth function in $h$. Hence, $\ang{(K\ast L)_x,\phi}$ is smooth in $x$, and we see that $K\ast L\in C_c^\infty(B,\EE'(G))$.

To see that $K\ast L$ is smooth away from $\GG^{(0)}\subset \GG$, write $K$ and $L$ as a sum
\[ K = K_\varepsilon + R_\varepsilon \;,\; L = L_\varepsilon + R'_\varepsilon ,\]
where $K_\varepsilon, L_\varepsilon$ have small propagation (i.e., they are supported in a small neighborhood of $\GG^{(0)} \subset \GG$), and  $R_\varepsilon, R'_\varepsilon$
are $C_c^\infty$ functions.
All one needs to verify is that $K_\varepsilon\ast R'_\varepsilon$ and $R_\varepsilon \ast L_\varepsilon$ are smooth functions on $\GG$.

\QED

A regular kernel $K$ on the groupoid $\GG$ is said to have an asymptotic expansion of degree $d$ near $B$,
\[ K \sim K^0+K^1+K^2+\cdots ,\]
if (1) each family $K^j$ is a regular kernel on $\GG$, (2) each $K_x^j$ is homogeneous of degree $-Q-d-j$ on $G$ and (3) the remainder term
\[ K - \sum_{j=0}^N K^j = R_N \] 
becomes more and more smooth as $N$ grows.
We denote the set of such kernels by $\Psi^d_H(\GG)$. 

\begin{proposition}
If $K\in \Psi^a_H(\GG)$ and $L\in \Psi^b_H(\GG)$ then $K\ast L\in \Psi^{a+b}_H(\GG)$.
In particular, $\Psi^0_H(\GG)$ is an algebra.
\end{proposition} 
(The proof is a simple adpatation of that found in [CGGP].)

If $K\in C_c^\infty(\GG)$ we have the left regular representation $\pi_x(K)$ of $K$ for each $x\in B$,
\[ \pi_x(K) \phi = K_x \ast \phi, \]
with $\phi \in C_c^\infty(G)$.
One regards $\pi_x(K)$ as a bounded operator on $L^2(G)$.
The $C^*$-norm on the convolution algebra of $\GG$ is defined,
for $K\in C_c^\infty(\GG)$, as
\[ \|K\|_{C^*(\GG)} = \sup_{x\in B} \, \|\pi_x(K) \| ,\]
and $C^*(\GG)$ is the norm closure of $C_c^\infty(\GG)$ with respect to this norm.
Now, the same definition applies to regular kernels of order zero.

\begin{proposition}
If $K\in \Psi_H^0(\GG)$ then $\pi_x(K)$ extends to a bounded operator on $L^2(G)$.
Moreover, the norm $\pi_x(K)$ is continuous in $x$.
\end{proposition}
Precisely as in the definition of the groupoid $C^*$-algebra we let
\[ \|K\|_{\Psi^0(\GG)} = \sup_{x\in B} \|\pi_x(K)\| ,\]
for $K\in \Psi^0_H(\GG)$.
We denote by $\PSI^0_H(\GG)$ the $C^*$-algebra that is the norm closure of $\Psi^0_H(\GG)$.

\section{Commutative diagrams}

In [Er1] we constructed a $K$-theory element $[\sigma_H(P)]\in K_0(C^*(T_HM)$,
and developed the tangent groupoid proof that shows that there is a natural map $K_0(C^*(T_HM))\to \ZZ$ that computes the index of $P$.
In section \ref{section:four}  of the present paper, using a short exact sequence, we showed that there is an element $[\sigma_H(P)]\in K_1(\overline{\calS}_H)$, and an index map $K_1(\overline{\calS}_H)\to \ZZ$.
In this section we connect the two approaches, by means of the order zero calculus on the tangent groupoid described in the previous section.

Let $\GG=\THM$ be the parabolic tangent groupoid of a compact contact manifold $M$.
Of course, the constructions described in the previous section can be carried out for $\GG=\THM$, since $\THM$ is locally identical to $\THG$.
We have a Rellich lemma.

\begin{proposition}
The norm closure of $\Psi^{-1}_H(\GG)$ in the left regular representation is $C^*(\GG)$.
In particular, $C^*(\GG)$ is a two-sided ideal in $\PSI^0_H(\GG)$.
\end{proposition}
{\bf Proof.}
The kernels $K\in \Psi^{-1}_H(\GG)$ are $L^1$ functions on $\GG$. 

\QED

The leading term in the asymptotic expansion of $K\in \Psi^0_H(\GG)$
can be regarded as a smooth family $K^0\in C^\infty([0,1],\calS_H)$.
We have a short exact sequence
\[ 0 \ra \Psi^{-1}_H(\GG) \ra \Psi^0_H(\GG) \ra C^\infty([0,1],\calS_H) \ra 0 .\]
Extending by continuity, we obtain an exact sequence of $C^*$-algebras,
\[ 0 \ra C^*(\GG) \ra \PSI^0_H(\GG) \ra C([0,1])\otimes \overline{\calS}_H \ra 0 .\]
Restriction to the $t=1$ boundary in the groupoid $\THM$ gives a commutative diagram
\[ \xymatrix{  0  \ar[r] & C^*(\GG) \ar[r]\ar[d] & \PSI^0_H(\GG) \ar[r]\ar[d] & C([0,1], \overline{\calS}_H) \ar[r]\ar[d] & 0 \\
               0  \ar[r] & \KK(L^2(M)) \ar[r] & \PSI^0_H(M) \ar[r] & \overline{\calS}_H \ar[r] & 0 }
\]
On the other hand, restriction to the $t=0$ boundary is a $K$-theory equivalence for each term in the sequence.
We obtain a commutative diagram in $K$-theory,
\[ \xymatrix{   K_1(\overline{\calS}_H) \ar[r]^{\partial}\ar[d]^{\cong}  & K_0(C^*(T_HM)) \ar[d]   \\
                K_1(\overline{\calS}_H) \ar[r]^{\partial} & K_0(\KK) = \ZZ. }
\]
Let's review the four maps that appear in this commutative square.
The horizontal maps $\partial$ are the boundary maps in $K$-theory. 
We already encountered the bottom one, and have seen that it sends $[\sigma_H(P)]\in K_1(\overline{\calS}_H)$ to the Fredholm index of $P$.
The horizontal map on the top of the diagram is the boundary map for the sequence
\[ 0 \ra C^*(T_HM) \ra \PSI^0_H(T_HM) \ra \overline{\calS}_H \ra 0 .\]
Here $\PSI^0(T_HM)$ is the closure in the left regular representation 
of the algebra $\Psi^0_H(T_HM)$ of regular kernels of order zero on the groupoid $T_HM$.

The class $[\sigma_H(P)]\in K_1(\calS_H)$ of a an order zero Fredholm operator is mapped to a class in $K_0(C^*(T_HM))$.
The vertical map on the left is the identity map.
Commutativity of the diagram therefore shows that the vertical map on the right
\[ K_0(C^*(T_HM)) \to \ZZ \]
sends the class $\partial(\sigma_H(P))\in K_0(C^*(T_HM))$ to the Fredholm index of $P$.

In [Er1], we described how the Connes-Thom isomorphism induces a natural isomorphism of $K$-theory classes
\[ K_0(C^*(G_m)) \cong K_0(C^*(T_mM)) \cong K^0(T^*_mM) ,\]
which, `glued' together, yield a natural isomorphism
\[ \Phi \;\colon\; K_0(C^*(T_HM)) \cong K^0(T^*M).\]
We showed that this isomorphism identifies the index map for Heisenberg symbols in $K_0(C^*(T_HM))$ with the topological index of Atiyah-Singer for elliptic symbols in $K^0(T^*M)$.
We summarize.

\begin{proposition}\label{prop:halfway}
Let $M$ be a closed contact manifold.
Let $P$ be an order zero pseudodifferential operator in the Heisenberg calculus on $M$ with invertible model operators $P_m$.

The boundary map in $K$-theory $\partial\colon K_1(\calS_H)\to K_0(C^*(T_HM))$ 
for the short exact sequence 
\[ 0 \ra C^*(T_HM) \ra \PSI^0_H(T_HM) \ra \overline{\calS}_H \ra 0 ,\]
composed with the Connes-Thom isomorphism $\Phi\colon K_0(C^*(T_HM))\cong K^0(T^*M)$
sends the invertible symbol $\sigma_H(P)\in \calS_H$ associated to $P$ to a $K$-theory class 
\[ \Phi(\partial(\sigma_H(P))) \in K^0(T^*M) .\]
Then
\[ {\rm Index}\, P = \int_{T^*M} {\rm Ch}(\Phi(\partial(\sigma_H (P))))\wedge {\rm Td}(M) .\]
\end{proposition}
Observe that, since $H\subseteq TM$ has a complex structure, it is a spin bundle, and so by the Thom isomorphism in $K$-theory
\[ K^0(T^*M) = K^0(H^*\oplus N^*) = K^0(N^*) .\]
Moreover, we assumed the existence of a global contact form, which means that $N^*\cong M\times \RR$.
Therefore $K^0(N^*) = K^1(M)$. 
To prove Theorem \ref{thm:main} we must show that the element $\Phi(\partial(\sigma_H(P)))\in K^0(T^*M)$ corresponds to the element $[a(P)]\in K^1(M)$.

\section{Proof of the main theorem, more diagrams}\label{section:nine}

We first prove Theorem \ref{thm:main} for a special class of operators called called {\em Hermite operators} (see [Ep], [EM]).
Let 
\[ {\cal I}\subseteq \Psi^0_H(M) \]
be the ideal of order zero operators $P$ on $M$ with $\pi_m(P_m^{op}) = 0$ for all $m\in M$.
The symbols $\sigma_H(P)$ of Hermite operators $P\in {\cal I}$ form an ideal $\sigma_H[{\cal I}]$ in $\overline{\calS}_H$.
Observe that, since $\pi_m(P_m)$ has the same principal Weyl symbol as $\pi_m(P_m^{op})$ (see section \ref{section:five}), it follows that $\pi_m(P_m)$ is a compact operator.
We therefore have a map
\[ \sigma_H[{\cal I}] \to \KK(V^{BF}) \;\colon\; \{P_m\} \mapsto \{\pi_m(P_m)\}. \]
The following lemma implies that this map is an isomorphism of $C^*$-algebras.

\begin{lemma}
For an order zero model operator $P_m$ we have
\[ \|P_m\| = \max \{\|\pi_m(P_m)\|,\|\pi_m(P_m^{op})\| \} .\]
\end{lemma}
{\bf Proof.}
Let $\pi_\tau, \tau\in \RR\setminus \{0\}$ denote the family of irreducible unitary representations of the Heisenberg group $G_m$, excluding the scalar representation.
Let $\pi_{+1}=\pi_m$ be the Bargmann-Fok representation, and $\pi_{-1}$ its conjugate representation.

Because $P_m$ is a distribution on $G_m$ homogeneous of order $-Q$, $\pi_\tau(P_m)$ is constant for $\tau>0$ as well as for $\tau<0$.
One verifies easily that $\|\pi_\tau(P_m^{op})\| = \|\pi_{-\tau}(P_m)\|$.
The lemma follows from the Plancherel formula for $G_m$.

\QED

The following theorem combined with Proposition \ref{prop:halfway} proves Theorem \ref{thm:main} for Fredholm operators $P$ with $\pi_m(P_m^{op}) = 1$.

\begin{theorem}\label{thm:hermite}
For a hypoelliptic  order zero operator $P$ with $P-1\in {\cal I}$, we have
\[ \Phi(\partial(\sigma_H(P))) = [a(P), V^{BF}] \in K^1(M)\cong K^0(T^*M).\]
\end{theorem}
{\bf Proof.}
Recall that the boundary map
\[ \partial \;\colon\; K_1(\overline{\calS}_H) \to K_0(C^*(T_HM)) \]
was associated to the sequence
\[ 0 \ra C^*(T_HM) \ra \PSI^0_H(T_HM) \ra \overline{\calS}_H \ra 0 .\]
Consider the two-sided ideal in $\PSI^0_H(T_HM)$ 
generated by regular order zero kernels $K_m$
for which $\pi_\tau(K_m)=0$ for all $\tau\le 0$.
This ideal can be naturally identified with $C((0,\infty],\KK(V^{BF}))$.
We have a commutative diagram
\[ \xymatrix{  0  \ar[r] & C_0((0,\infty),\KK(V^{BF})) \ar[r]\ar[d] & C_0((0,\infty],\KK(V^{BF})) \ar[r]\ar[d] & \KK(V^{BF}) \ar[r]\ar[d] & 0 \\
               0  \ar[r] & C^*(T_HM) \ar[r] & \PSI^0_H(T_HM) \ar[r] & \overline{\calS}_H \ar[r] & 0 }
\]
which induces a commutative diagram in $K$-theory,
\[ \xymatrix{ K_1(\KK(V^{BF})) \ar[r]^{\partial}\ar[d] & K_0(C_0((0,\infty),\KK(V^{BF}))) \ar[d] \\
							K_1(\overline{\calS}_H) \ar[r]^{\partial} & K_0(C^*(T_HM)) . }
\]
The boundary map on the top row is just the standard isomorphism $K_1(A)\cong K_0(C_0(\RR,A))$. 

We try to understand what the Connes-Thom isomorphism does to the $K$-theory of the subalgebra
\[ C_0((0,\infty),\KK(V^{BF})) \subseteq C^*(T_HM).\]
As one easily verifies, the deformation that induces the Connes-Thom isomorphism (associated to the adiabatic groupoid of $T_HM$, see [Er1]) gives rise to a commutative diagram
\[ \xymatrix{ K_0(C_0((0,\infty),\KK(V^{BF}))) \ar[r]\ar[d]  & K_0(C_0((0,\infty),C_0(H^*))) \ar[d] \\ 
							K_0(C^*(T_HM)) \ar[r] & K_0(C_0(H^*\oplus N^*)). }
\] 
The map on the top row is induced fiberwise by the Weyl quantization map from $C_0(H^*_m)$ to $\KK(V^{BF}_m))$, which induces the Bott isomorphism in $K$-theory.
The full map thus corresponds precisely to the ordinary Thom isomorphism for the complex bundle $H^*\times \RR$ over the base space $M\times \RR$
\[ K^0(M\times\RR) \cong K^0(H^*\times\RR) .\] 
Notice that the vertical map on the right is induced by the inclusion of $H^*\times(0,\infty)$
into $H^*\times \RR \cong H^*\oplus N^*$. 
This inclusion is homotopic to the identity, and induces an isomorphism in $K$-theory.
 
Combining both commutative squares, we obtain
\[
 \xymatrix{  K_1(\KK(V^{BF})) \ar[r]^{\cong}\ar[d]_{\subset} & K_0(\KK(V^{BF})\otimes C_0(\RR)) \ar[r]^{\cong}\ar[d]^{\subset} & K_0(C_0(H^*)\otimes C_0(\RR)) \ar[d]^{\cong} \\
             K_1(\overline{\calS}_H) \ar[r]_{\partial} & K_0(C^*(T_HM)) \ar[r]_{\Phi} & K^0(T^*M) }.
\]
Now, if $P-1\in {\cal I}$, we have simply $a(P)_m=\pi_m(P_m)$, and the theorem follows  from commutativity of the diagram.

\QED

The general case reduces to that of Hermite operators by means of a symmetry argument.

\begin{proposition}
Let  $\sigma\in K_1(\overline{\calS}_H)$.
Then $\Phi(\partial(\sigma^{op})) = - \Phi(\partial(\sigma))$.
\end{proposition} 
{\bf Proof.}
The anti-automorphisms `op' of the osculating groups was defined as
\[ {\rm op} \;\colon\; H_m\oplus N_m \to H_m\oplus N_m \;\colon\; (h,t)\mapsto (h,-t) .\]
These maps induce anti-automorphisms of the algebra $\PSI^0_H(T_HM)$,
its  ideal $C^*(T_HM)$, and the quotient $\overline{\calS}_H$.
Under the Connes-Thom isomorphism
\[ \Phi\;\colon\; K_0(C^*(T_HM)) \to K^0(H^*\oplus N^*) \]
the map `op' corresponds to the obvious map in $K^0(H^*\oplus N^*)$ induced by $(h,t)\mapsto (h,-t)$.
Therefore, in $K^0(H^*\oplus N^*)\cong K^0(N^*)$ the map `op' reverses the sign of the $K$-theory class.

\QED

It follows that if $\sigma^{op}=\sigma$ in $K_1(\overline{\calS}_H)$ 
then $\Phi(\partial(\sigma))=0$.
Theorem \ref{thm:main} now follows immediately.

\section*{References}

\def\item{\vskip2.75pt
plus1.375pt minus.6875pt\noindent\hangindent1em} \hbadness2500
\tolerance 2500
\markboth{References}{References}

\item{[BG]} R.\ Beals and P.\ Greiner,
{\sl Calculus on Heisenberg Manifolds},
Annals of Mathematics Studies (119), Princeton, 1988.

\item{[Bo]} L.\ Boutet de Monvel, 
{\sl On the index of Toeplitz operators of several complex variables},
Invent. Math. 50 (1979), 249--272.

\item{[CGGP]}, M.\ Christ, D.\ Geller, P.\ Glowacki, L.\ Polin,
{\sl Pseudodifferential operators on Groups with dilations},
Duke Math.\ J\. 68 (1992), 31--65.

\item{[EM]} C.\ Epstein and R.\ Melrose,
{\sl Contact degree and the index of Fourier integral operators},
Math. Res. Lett. 5 (1998), no.3, 363--381.

\item{[EM2]} C.\ Epstein and R.\ Melrose,
{\sl The Heisenberg algebra, index theory and homology},
preprint, 2003.

\item{[Ep]} C.\ Epstein,
{\sl Lectures on Indices and Relative Indices on Contact and CR-manifolds},
Woods Hole Mathematics: Perspectives in Mathematics and Physics,
World Scientific, 2004.

\item{[Er1]} E. van Erp,
{\sl The Atiyah-Singer Formula for Subelliptic Operators on a Contact Manifold, part I},
in preparation. 

\item{[FS1]} G.\ B.\ Folland and E.\ M.\ Stein,
{\sl Estimates for the $\bar{\partial}_b$ complex and analysis on the Heisenberg group},
Comm.\ Pure and Appl.\ Math.\, vol XXVII (1974), 429--522.

\item{[Ta1]} M.\ E.\ Taylor,
{\sl Noncommutative microlocal analysis, part I},
Mem.\ Amer.\ Math.\ Soc.\, vol. 313, AMS, 1984.

\end{document}